\documentclass[12pt,reqno]{amsart}

\usepackage[foot]{amsaddr}
\usepackage{amsfonts}
\usepackage{amsmath}
\usepackage{amssymb}
\usepackage{amsthm}
\usepackage[french,english]{babel}
\usepackage{bm}
\usepackage{booktabs}
\usepackage{colortbl}
\usepackage{diagbox}
\usepackage{epsfig}
\usepackage{epstopdf}
\usepackage{enumerate}
\usepackage{enumitem}
\usepackage[T1]{fontenc}
\usepackage{framed}
\usepackage{fullpage}
\usepackage[left=2.50cm, right=2.50cm, top=2.50cm, bottom=2.50cm]{geometry}
\usepackage{graphicx}
\usepackage{hyperref}
\usepackage{lastpage}
\usepackage{latexsym}
\usepackage{mathtools}
\usepackage{mathrsfs}
\usepackage{multirow}
\usepackage{pstricks}
\usepackage{subfig}
\usepackage{tikz}
\usepackage{tikz-cd}
\usepackage{tikz-qtree}
\usepackage{todonotes}
\usepackage{url}
\usepackage{verbatim}
\usepackage{xcolor}

\usetikzlibrary{arrows}
\usetikzlibrary{automata}
\usetikzlibrary{cd}
\usetikzlibrary{matrix}
\usetikzlibrary{positioning}
\usetikzlibrary{patterns}

\usepackage{cleveref}
\hypersetup{hypertex=true,colorlinks=true, 
linkcolor=red,
menucolor=cyan,
citecolor=cyan,
filecolor=cyan,
urlcolor=blue
}

\def\centerarc[#1](#2)(#3:#4:#5)
{ \draw[#1] ($(#2)+({#5*cos(#3)},{#5*sin(#3)})$) arc (#3:#4:#5); }

\makeatletter
\renewcommand{\email}[2][]{%
  \ifx\emails\@empty\relax\else{\g@addto@macro\emails{,\space}}\fi%
  \@ifnotempty{#1}{\g@addto@macro\emails{\textrm{(#1)}\space}}%
  \g@addto@macro\emails{#2}%
}
\makeatother

\makeatletter

\@namedef{subjclassname@2010}{

  \textup{2020} Mathematics Subject Classification}

\makeatother

\newtheorem{thm}{Theorem}[section]

\newtheorem*{thm*}{Theorem}

\newtheorem{prop}[thm]{Proposition}

\theoremstyle{definition}

\newtheorem*{rmk}{Remark}

\numberwithin{equation}{section}

\newcommand{\newabstract}[1]{%
  \par\bigskip
  \csname otherlanguage*\endcsname{#1}%
  \csname captions#1\endcsname
  \item[\hskip\labelsep\scshape\abstractname.]
}
\makeatletter
\tikzset{
        hatch distance/.store in=\hatchdistance,
        hatch distance=5pt,
        hatch thickness/.store in=\hatchthickness,
        hatch thickness=5pt
        }
\pgfdeclarepatternformonly[\hatchdistance,\hatchthickness]{north east hatch}
    {\pgfqpoint{-1pt}{-1pt}}
    {\pgfqpoint{\hatchdistance}{\hatchdistance}}
    {\pgfpoint{\hatchdistance-1pt}{\hatchdistance-1pt}}%
    {
        \pgfsetcolor{\tikz@pattern@color}
        \pgfsetlinewidth{\hatchthickness}
        \pgfpathmoveto{\pgfqpoint{0pt}{0pt}}
        \pgfpathlineto{\pgfqpoint{\hatchdistance}{\hatchdistance}}
        \pgfusepath{stroke}
    }
\makeatother

\begin{document}
\title[A Note on $\zeta(5)$]{A Note on a Recent Attempt to Prove the Irrationality of $\zeta(5)$}

\author{Keyu Chen}
\author{Wei He}
\author{Yixin He}
\author{Yuxiang Huang}
\author{Yanyang Li}
\author{Quanyu Tang}
\author{Lei Wu}
\author{Shenhao Xu}
\author{Shuo Yang}
\author{Zijun Yu}
\address[Keyu Chen]{School of Mechanical Engineering, Sichuan University, Chengdu 610064, P. R. China}
\address[Wei He]{School of Mathematical Sciences, Chongqing Normal University, Chongqing 401331, P. R. China}
\address[Yixin He]{College of Mathematics and Statistics, Northwest Normal University, Lanzhou 730000, P. R. China}
\address[Yuxiang Huang]{School of Mathematics and Statistics, Lanzhou University, Lanzhou 730000, P. R. China}
\address[Yanyang Li]{School of Mathematics, Southeast University, Nanjing 210000, P. R. China}
\address[Quanyu Tang, Zijun Yu]{School of Mathematics and Statistics, Xi'an Jiaotong University, Xi'an 710049, P. R. China}
\address[Lei Wu]{School of Mathematics, Sichuan University, Chengdu 610064, P. R. China}
\address[Shenhao Xu]{Technical College for the Deaf, Tianjin University of Technology, Tianjin 300384, P. R. China}
\address[Shuo Yang]{School of Mathematics and Information Sciences, Yantai University, Yantai 264000, P. R. China}

\email{2770891488@qq.com}
\email{weih5907@gmail.com}
\email{hyx717math@163.com}
\email{hyuxiang2023@lzu.edu.cn}
\email{liyanyang1219@gmail.com}
\email{tang\_quanyu@163.com}
\email{leiwu0829@163.com}
\email{empetmb@163.com}
\email{shuo96033@gmail.com}
\email{yuzijun0830@outlook.com}

\date{}

\begin{abstract}
Recently, Shekhar Suman \cite{Suman2024} made an attempt to prove the irrationality of $\zeta(5)$. Unfortunately, the proof is not correct. In this note, we discuss the fallacy in the proof.
\end{abstract}

    	\maketitle

\section{Introduction}

The study of irrationality and transcendence of numbers is a classical topic in transcendental number theory. Recall the Riemann zeta values
\[
\zeta(n) = \sum_{k \geq 1} \frac{1}{k^n} \quad \text{for } n \geq 2.
\] Euler proved that
\[
\zeta(2) = \frac{\pi^2}{6}
\]
and more generally
\[
\zeta(2n) = -\frac{B_{2n} (2\pi i)^{2n}}{2 (2n)!} \quad \text{for } n \geq 1,
\]
where \( B_m \) is the \( m \)-th Bernoulli number. It is conjectured that the values of the Riemann zeta function at all positive integers are irrational.  

In 1882, Lindemann proved that the number \(\pi\) is transcendental (see, for example, \cite[Theorem 1.2]{Waldschmidt2013}). As a consequence, the even values \(\zeta(2n)\) are irrational. In 1979, Ap\'{e}ry \cite{Apery1979} famously demonstrated that the number \(\zeta(3)\) is irrational. Subsequently, Beukers \cite{Beukers1979} provided an elementary proof of \(\zeta(3)\) based on integral representations. In 2001, Keith Ball and Tanguy Rivoal \cite[Th\'{e}or\`{e}me
 1]{BallRivoal2001} applied Nesterenko’s criterion to prove that the \(\mathbb{Q}\)-vector space spanned by odd zeta values is infinite-dimensional:
\[
\operatorname{dim}_{\mathbb{Q}} \langle 1, \zeta(3), \zeta(5), \ldots, \zeta(2n+1) \rangle_{\mathbb{Q}} \geq  \frac{1}{3}\log(2n+1).
\] In the same year, Zudilin \cite{Zudilin2001} demonstrated that at least one of the four numbers
\[
\zeta(5), \zeta(7), \zeta(9), \zeta(11)
\]
is irrational. Incidentally, the proof of the irrationality of \(\zeta(3)\) can also be found in Zudilin's recent book \cite[Section 7]{Zudilin2023}.

To this day, it remains unknown whether \(\zeta(5) \notin \mathbb{Q}\). It is also not known whether \(\zeta(3) \notin \pi^3 \mathbb{Q}\).

Recently, Shekhar Suman \cite{Suman2024} claimed to have proved the irrationality of \(\zeta(5)\). Unfortunately, the proposed proof is incorrect.  

This note is organized as follows. In Section \ref{s2}, we identify the logical flaws in the proof presented in \cite{Suman2024}. In Section \ref{s3}, we outline a standard approach to proving the irrationality of a number and highlight that the irrationality of \(\zeta(2m+1)\) for general odd integers \(2m+1\) remains an open and intriguing area of research. In Section \ref{s4}, we discuss recent advancements, including the groundbreaking work of Calegari, Dimitrov, and Tang \cite{calegari2024linearindependence1zeta2}, which introduces a powerful new method for proving irrationality and offers promising directions for future research.

\section{Logical Fallacies in Suman's Proof}\label{s2}

In this section, we analyze the flaws in Shekhar Suman's proof of \cite[Theorem 1]{Suman2024}. While the preliminary result, \cite[Lemma 1]{Suman2024}, is correct and well-established, the proof of \cite[Theorem 1]{Suman2024} contains a critical logical error.

Suman claims that the linear Diophantine equation \cite[Eq. (47)]{Suman2024} has no integer solution, which he uses to derive the irrationality of \(\zeta(5)\). However, we show that \cite[Eq. (47)]{Suman2024} indeed admits integer solutions, invalidating his conclusion.

Suman argues as follows:
\begin{quote}
Since \((d_n, 2d_n) = d_n\), the above linear Diophantine equation has an integral solution if and only if \(d_n \mid k_i b\). So we have
\[
d_n a - 2d_n b = -k_i b \quad \text{where } d_n \mid k_i b, \quad 1 \leq k_i \leq d_n - 1, \quad n \geq b.
\]
He then claims that the following cases are impossible:
\[
d_n a - 2d_n b = -k_i b \quad \text{where } d_n \mid k_i b, \quad 0 \leq k_i \leq d_n, \quad n \geq 1.
\]
\end{quote}

He employs induction on \(n\) to establish this claim. However, a closer inspection reveals a critical flaw in his argument at the base case \(n=1\).

For \(n=1\), the equation becomes:
\[
a - 2b = -k_i b, \quad \text{where } 0 \leq k_i \leq 1 \text{ and } 1 \mid k_i b.
\]
Solving this, Suman correctly concludes that the possible solutions are:
\[
a - 2b = 0 \quad \text{or} \quad a - 2b = -b.
\]
This implies that \(a = 2b\) or \(a = b\), both of which are valid integer solutions. However, instead of acknowledging the existence of these integer solutions as a refutation of his claim, Suman proceeds to argue that these solutions lead to:
\[
\zeta(5) = \frac{a}{b} = 2 \quad \text{or} \quad \zeta(5) = 1,
\]
which he deems "absurd" because it is well-known that \(\zeta(5)\) is not an integer.

Here lies the critical error: the claim that equation \cite[Eq. (48)]{Suman2024} has no integer solutions is logically independent of whether \(\zeta(5)\) is an integer or not. The well-established fact that \(\zeta(5)\) is not an integer cannot be used to prove that \cite[Eq. (48)]{Suman2024} has no integer solutions. The existence of integer solutions for \cite[Eq. (48)]{Suman2024} at \(n=1\) directly invalidates the induction base case and, consequently, the entire proof.

To summarize, the error in \cite[Theorem 1]{Suman2024} lies in conflating two unrelated facts: the solvability of the Diophantine equation \cite[Eq. (48)]{Suman2024} and the irrationality of \(\zeta(5)\). The former is a purely algebraic property of the equation, while the latter is a number-theoretic property of the zeta function. By assuming \(\zeta(5)\) is not an integer to argue against the solvability of \cite[Eq. (48)]{Suman2024}, Suman undermines the logical foundation of his proof. As a result, Theorem 1 in \cite{Suman2024} is incorrect.

In addition to the errors in the proof of \cite[Theorem 1]{Suman2024}, similar logical fallacies can be identified in the proof of \cite[Theorem 2]{Suman2024}, which aims to establish the irrationality of \(\zeta(2m+1)\) for integers \(m \geq 2\). While the specific details of the argument in \cite[Theorem 2]{Suman2024} differ from \cite[Theorem 1]{Suman2024}, the underlying issue remains the same: the reliance on assumptions or conclusions that are independent of the solvability of the corresponding Diophantine equations.

\section{A Standard Approach to Proving the Irrationality of a Number}\label{s3}

In this section, we review a standard method for proving the irrationality of a number \(\alpha\) (the criterion for irrationality), and show that Suman's method does not meet the criterion.

\begin{prop}\label{p1}
Suppose that we can construct sequences of pairs of rational numbers \(a_n, b_n\) with the following properties:
\begin{itemize}
    \item (1) There is a small number \(0 < \varepsilon < 1\) such that
    $$
    0 < \left| a_n \alpha - b_n \right| < \varepsilon^n
    $$
    for all sufficiently large \(n\).
    \item (2) Let \(d_n \in \mathbb{N}\) be the common denominator of \(a_n, b_n\):
    $$
    d_n a_n \in \mathbb{Z}, \quad d_n b_n \in \mathbb{Z}.
    $$
    Assume that \(d_n < D^n\) for some \(D \in \mathbb{R}\).
    \item (3) \(D\) is not too big:
    $$
    D \varepsilon < 1.
    $$
\end{itemize}
Then \(\alpha\) is irrational.
\end{prop}

\begin{proof} (by contradiction) Suppose that \(\alpha\) is rational, \(\alpha = \frac{p}{q}\) where \(p, q \in \mathbb{Z}, q > 0\). Assumption (1) then becomes
$$
0 < \left| a_n \frac{p}{q} - b_n \right| < \varepsilon^n \quad \text{for large } n
.$$ By multiplying through by \(q\) and \(d_n\), we obtain
$$
0 < \left| d_n a_n p - d_n b_n q \right| < q d_n \varepsilon^n < q D^n \varepsilon^n
.$$ Since by assumption (3), \(D \varepsilon < 1\), the right-hand side tends to zero. Thus, we can find a large \(n\) such that
$$
0 < \left| \underbrace{(d_n a_n)}_{\in \mathbb{Z}} p - \underbrace{(d_n b_n)}_{\in \mathbb{Z}} q \right| < 1
.$$ But by (2), this is an integer between 0 and 1, which is a contradiction.
\end{proof}

\begin{rmk}
    Besides Beukers \cite{Beukers1979}, other interesting applications of Proposition \ref{p1} can be found in, for example, Huylebrouck \cite{Huylebrouck2001}.
\end{rmk}

Numerical verification shows that the \(I_n\) in \cite[Lemma 1]{Suman2024} does not satisfy the condition (3) of Proposition \ref{p1}. Therefore, \(I_n\) cannot be used to prove the irrationality of \(\zeta(5)\) via Proposition \ref{p1}. In fact, mathematicians have been searching for a suitable integral representation of the form \(a_n + b_n \zeta(5)\) that converges to 0, but such a representation has not yet been found. Therefore, the irrationality of \(\zeta(2m+1)\) for \(m \geq 2\) still remains an open problem worth exploring. Finally, finding a good integral representation for \(\zeta(5)\) or even for \(\zeta(2m+1)\) that can be used to approximate them continues to be a worthwhile pursuit.

\section{Future Perspectives}\label{s4}

Recently, Calegari, Dimitrov, and Tang \cite{calegari2024linearindependence1zeta2} made significant progress in the study of irrationality and $\mathbb{Q}$-linear independence of special values in number theory. Their work established several groundbreaking results, including the irrationality of the classical Dirichlet $L$-value
\[
L(2, \chi_{-3}) = \sum_{n=0}^\infty \left(\frac{1}{(3n+1)^2} - \frac{1}{(3n+2)^2}\right),
\]
and demonstrated its linear independence with $1$ and $\pi^2$ over $\mathbb{Q}$. This result, formalized as \cite[Theorem A]{calegari2024linearindependence1zeta2} in their paper, represents a significant advancement in our understanding of special values of $L$-functions. Furthermore, their techniques generalized Zagier's constructions and provided a new arithmetic holonomy bound framework.

Their work also yielded \cite[Corollary B]{calegari2024linearindependence1zeta2}, which proves the irrationality of numerous related series. For instance, they showed that certain linear combinations of $L(2, \chi_{-3})$, $\pi^2$ are irrational. These results expand the set of known irrational numbers derived from $L$-functions and highlight the potential of their methods for exploring similar problems.

Moreover, \cite[Theorem C]{calegari2024linearindependence1zeta2} extends their analysis to logarithmic values, establishing the irrationality of \( \log\left(1 + \frac{1}{m}\right) \log\left(1 + \frac{1}{n}\right)\) under the condition \( |\frac{m}{n}-1|<\frac{1}{10^6}\) and proving the $\mathbb{Q}$-linear independence of the following four expressions:
\[
1, \quad
\log\left(1 + \frac{1}{m}\right), \quad \log\left(1 + \frac{1}{n}\right), \quad
\log\left(1 + \frac{1}{m}\right) \log\left(1 + \frac{1}{n}\right).
\]

These findings represent a major breakthrough in the study of special values of zeta and $L$-functions, with implications for broader problems in number theory, transcendence theory, and the study of modular forms and related automorphic objects. The techniques introduced by Calegari, Dimitrov, and Tang provide a robust foundation for future exploration, including the study of other special values and related constants in mathematics.

\end{document}